\documentclass[12pt]{article}
\usepackage[latin1]{inputenc}
\usepackage{fullpage}
\usepackage{amsfonts}
\usepackage{latexsym,amssymb}
\usepackage{amsmath}
\usepackage{subfigure}
\usepackage{graphicx}        % standard LaTeX graphics tool
\usepackage{epsf}

\title{Environmental Noise Variability\\
in Population Dynamics Matrix Models}
\author{Michel \textsc{De Lara}\\
Universit\'e Paris-Est, Cermics,\\
6-8 avenue Blaise Pascal, 77455 Marne la Vall\'ee Cedex 2, France\\
delara( )cermics.enpc.fr}

\newtheorem{proposition}{Proposition}

\def\Figdir{./}

\newcommand{\EE}{\mathbb{E}}
\newcommand{\PP}{\mathbb{P}}
\newcommand{\RR}{\mathbb{R}}
\newcommand{\MM}{\mathbb{M}}

\newcommand{\greater}{\succeq} % 
\newcommand{\lesser}{\preceq}

\def\defegal{:=}
\newcommand{\vari}{\mathop{\mathsf{var}}}
\newcommand{\CV}{\mathop{\mathsf{CV}}}

\def\abundance{N}

\def\abundance{n}
\def\size{N}
\def\horizon{T}
\def\matrix{A}

\def\alert#1{#1}
\def\pause{}
\def\mtext#1{\,\mbox{#1}\,} %text in maths

\newenvironment{proof}{\small{\bf Proof.}}{\hfill$\Box$\normalsize
\bigskip}

\begin{document}
\maketitle

\begin{abstract}
The impact of  environmental variability on population size growth rate
in dynamic models is a recurrent issue in the theoretical
ecology literature.
In the scalar case, R.~Lande pointed out that results are ambiguous
depending on whether the noise is added at arithmetic or logarithmic
scale, while the matrix case has been investigated by S.~Tuljapurkar.
Our contribution consists first in introducing another notion of
variability than the widely used variance or coefficient of
variation, namely the so-called convex orders.
Second, in population dynamics matrix models, we focus on how matrix
components depend functionaly on uncertain environmental factors.
In the log-convex case, we show that, in a sense, 
environmental variability increases both 
mean population size and mean log-population size 
and makes them more variable.
Our main result is that
specific analytical dependence coupled with appropriate notion
of  variability lead to wide generic results, valid for all times and
not only asymptotically, and requiring no assumptions of
stationarity, of normality, of independency, etc. 
Though the approach is different, our conclusions are consistent
with previous results in the literature.
However, they make it clear that the analytical dependence on 
environmental factors cannot be overlooked when trying to tackle the
influence of variability.
\end{abstract}

\begin{quote}

\end{quote}

\pagebreak

\tableofcontents

\section{Influence of environmental noise on population size}

We recall here different observations and results in the theoretical
ecology literature which point out the ambiguous role of environmental noise
on population size in matrix population models, 
according to whether the 
noise is added at arithmetic or logarithmic scale.

\subsection{Lande's comments on additive noise at
arithmetic or logarithmic scale}

R.~Lande in \cite{Lande-Engen-Saether:2003} comments the influence of 
environmental noise on population size according to whether the 
noise is added at arithmetic or logarithmic scale.
The evolution of \alert{population size} $\size(t) $ in absence of
density-dependent effect may be described 
\begin{itemize}
 \item either on  \pause \alert{arithmetic scale} with 
\alert{multiplicative growth rate} $\lambda(t)$ and
 \pause \alert{dynamic}\\ $\size(t+1) = \lambda(t) \size(t) $,
\item or on \alert{logarithmic scale} with
\alert{growth rate on the log scale} 
$r(t)=\log \big(\lambda(t) \big)$ and
\alert{dynamic on the log scale} 
$\log \size(t+1) = r(t) + \log \size(t) $.
\end{itemize}
On the one hand, 
adding environmental noise to multiplicative growth rate as in
\(
\alert{\lambda(t) = \overline{\lambda} + \epsilon(t) }
\),
where the noise is zero-mean (\( \EE[\epsilon(t)]=0 \)),
gives the following mean of growth rate on the log scale
\[
\alert{ \overline{r} = \EE[ \log \big(\lambda(t) \big) ] 
\approx \log \overline{\lambda}   - \sigma_r^2 } \; .
\]
``Thus, demographic and \alert{environmental stochasticity reduce the mean
growth rate of a population on the logarithmic scale}, compared with
that in the (constant) average environment'' 
\cite{Lande-Engen-Saether:2003}.

On the other hand, 
adding environmental noise to growth rate on the log scale as in 
\(
\alert{r(t) = \overline{r} + \epsilon(t) }
\)
gives, in case $\epsilon(t)$ follows a Normal distribution 
${\cal N}(\overline{\epsilon}, \sigma_{\epsilon}^2)$,  
the following mean of growth rate on the arithmetic scale
\[
\alert{ \overline{\lambda} = \exp \big( 
\overline{r} + \overline{\epsilon} + \frac{\sigma_{\epsilon}^2}{2} 
\big) }  \; .
\]
Thus, Lande concludes that,
``with the mean environmental effect equal to zero, 
$ \overline{\epsilon}= 0 $, then it would be found that 
\alert{environmental  stochasticity increases the mean multiplicative
  growth rate}, $ \overline{\lambda} $''.

\subsection{Tuljapurkar's asympotic approximation}

S.~Tuljapurkar considers a stationary sequence of random matrices
$\matrix_0$, $\matrix_1$, \ldots yielding population vector
$\abundance(t)= \matrix_{t-1} \cdots \matrix_0 \abundance(0) $ and 
population size $\size(t) = \| \matrix_{t-1} \cdots \matrix_0 \abundance(0) \| $.
Under general conditions (see \cite{Tuljapurkar:1990,Caswell:2001}),
there exists a deterministic \emph{stochastic growth rate}
\( \lambda_s \) defined by
\[
\log \lambda_s = 
\lim_{t \to +\infty} \frac{1}{t} \log \size(t) 
 = \lim_{t \to +\infty}
\frac{1}{t} \log \| \matrix_{t-1} \cdots \matrix_0 \abundance(0) \| \; .
\]
% What is more, the distribution of the population size $\size(t)$
% may be shown to be asymptotically lognormal
% (see \cite{Furstenberg and Kesten (1960), Orzack-Tuljapurkar:1980})
% {
% \[
%  \log \size(t) \to_{t \to +\infty} {\cal N}\!\textrm{ormal} 
% (t \log \lambda_s , t \sigma^2 ) \; .
% \]
% }
Denoting by $\lambda_1$ the largest eigenvalue of the 
{average matrix} $\overline{\matrix}$, Tuljapurkar obtains the approximation
{
\[
\log \lambda_s \approx \log \lambda_1 \pause - \frac{\tau^2}{2\lambda_1^2}
\pause + \frac{\theta}{\lambda_1^2}
\]
}
where $\tau^2$ is proportional to the variance 
{$\EE[ (\matrix_t-\overline{\matrix}) \otimes (\matrix_t-\overline{\matrix}) ]$}
(and $\theta$ is related to autocorrelation).
In this case, environmental stochasticity reduces the mean
growth rate of the population.

\subsection{A quest for generic results}

The two above cases show that  environmental noise has 
an ambiguous impact on population size in matrix population models.
Our main objective is contributing to clarify this impact with generic
mathematical results. 
For this, we shall first introduce in Sect.~\ref{sec:Convex_orders} 
a tool to measure variability, distinct from the widely used variance or 
coefficient of variation, and known as \emph{convex partial orders}.
Then, in Sect.~\ref{sec:Generic_results}, 
we shall provide generic results on environmental noise variability
in population dynamics matrix models.
We conclude in Sect.~\ref{sec:Conclusion} by pointing out 
proximities and differences between our approach and those presented
in Sect.~\ref{sec:Convex_orders}.

\section{Convex orders as tools for measuring variability}
\label{sec:Convex_orders}

To a (square integrable) random variable $X$, one can attach the 
variance $\vari(X)$. This latter scalar measures ``variability'',
and any pair of random variables $X$ and $Y$ may be compared,
with $X$ being more variable than $Y$ if $\vari(X) \geq \vari(Y)$.
The variance thus defines a \emph{total order}.

Other orders are interesting for comparing pairs of random variables.
However, they are generally not total: not all pairs may be ranked.
Related to this is the fact that no single scalar, such as variance, 
may be attached to a random variable to measure its variability.
In this vein, we shall present the so-called \emph{increasing convex} 
and \emph{convex stochastic orders}.
Such orders can only rank random variables for which 
the primitives of their respective repartition functions never cross.

We think that these orders and many others referenced in the two main
books \cite{Muller-Stoyan:2002,Shaked-Shanthikumar:2007} 
may be useful in the ecological modelling scientific community. 
Of course, for this, practical tests must be developed to
compare empirical data as to their variability. 
This is not the object of this paper.

All random variables are defined on a probability space with probability 
$\PP$.
To a random variable $X$, we shall attach its 
(right-continuous) \emph{repartition function} $F(x)=\PP(X \leq x)$.
% The right-continuous version of its inverse is
% $F^{-1}(p)=\sup \{ x \mid F(x) \leq p \}$. 
We shall always consider random variables with finite means,
with generic notation $X$ and $Y$, 
and $F$ and $G$ for their respective repartition functions.

\subsection{Increasing convex order}

The increasing convex order compares random variables according both to
their ``location'' and to their
``variability'' or ``spread'' \cite{Shaked-Shanthikumar:2007}.
% \cite[p.16]{Muller-Stoyan:2002}
We say that $X$ \emph{is less than} $Y$ 
\emph{in increasing convex order}, denoted by
\begin{equation*}
X \lesser_{icx} Y \; ,
\end{equation*}
if and only if one of the following equivalent conditions holds true
\begin{itemize}
\item  the primitive of the repartition function of $X$ is always 
below that of $Y$: \\ 
$ \int_{-\infty}^{c} F(x) dx \leq \int_{-\infty}^{c} G(x) dx$,
for all $c \in \RR$,
\item $\EE ( \varphi (X ) )  \leq \EE ( \varphi (Y) )$ 
for all increasing and convex function $\varphi$. 
\end{itemize}
Roughly speaking, $Y$ is more likely to take on extreme values than $X$.
In a sense, $X$ is both ``smaller'' and ``less variable'' than $Y$
\cite{Shaked-Shanthikumar:2007}. 
We have the important property that,
when $X \lesser_{icx} Y $, the means are ordered too:
$\EE(X) \leq \EE(Y)$.
However, nothing can be said of the variances.
To compare variances, we need a stronger (more demanding) order. 

\subsection{Convex order}

\begin{figure}\centering
\begin{tabular}[p]{cc}
\subfigure[Repartition functions]{%
\epsfxsize=7cm\epsfbox{\Figdir 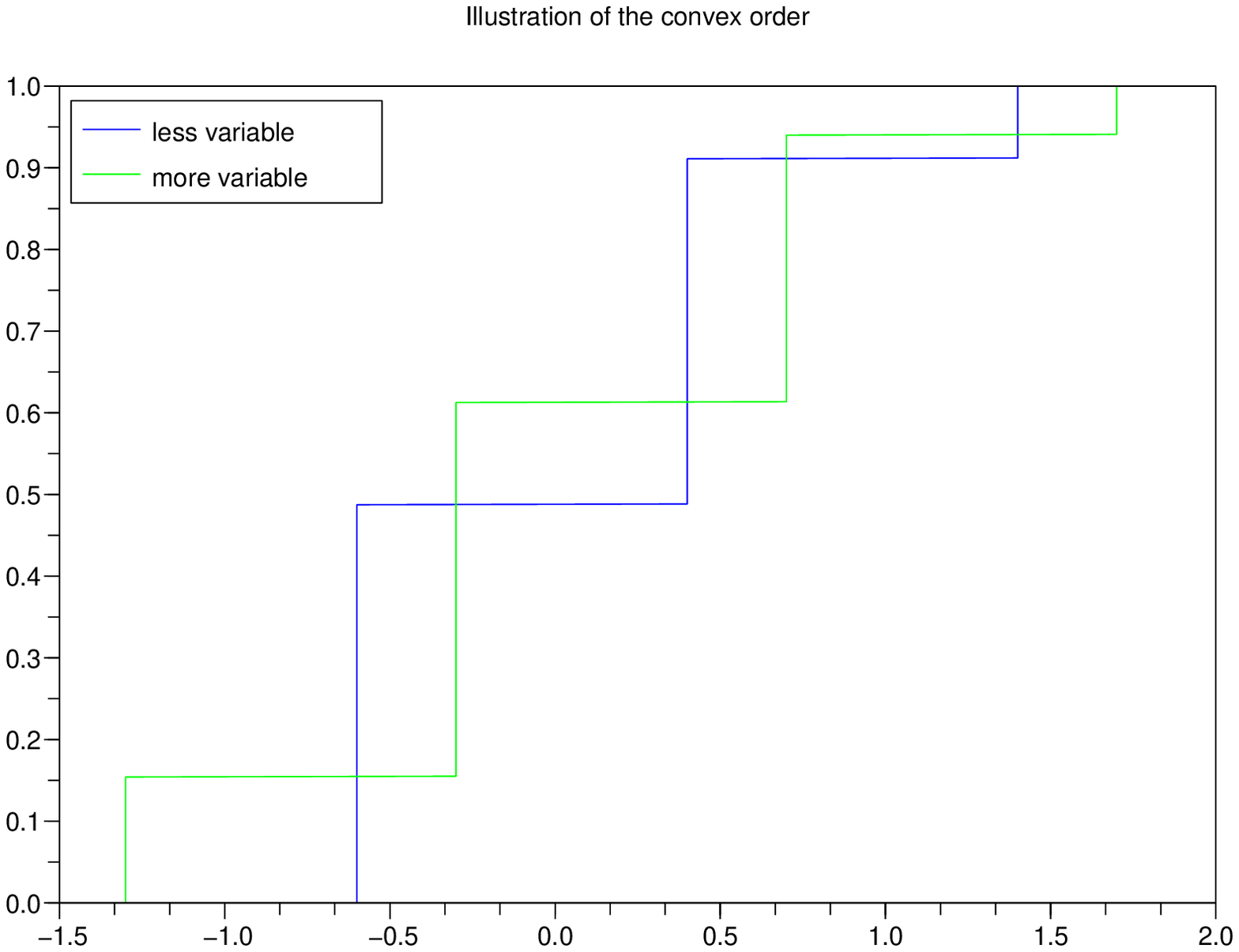}} &
\subfigure[Primitives]{%
\epsfxsize=7cm\epsfbox{\Figdir 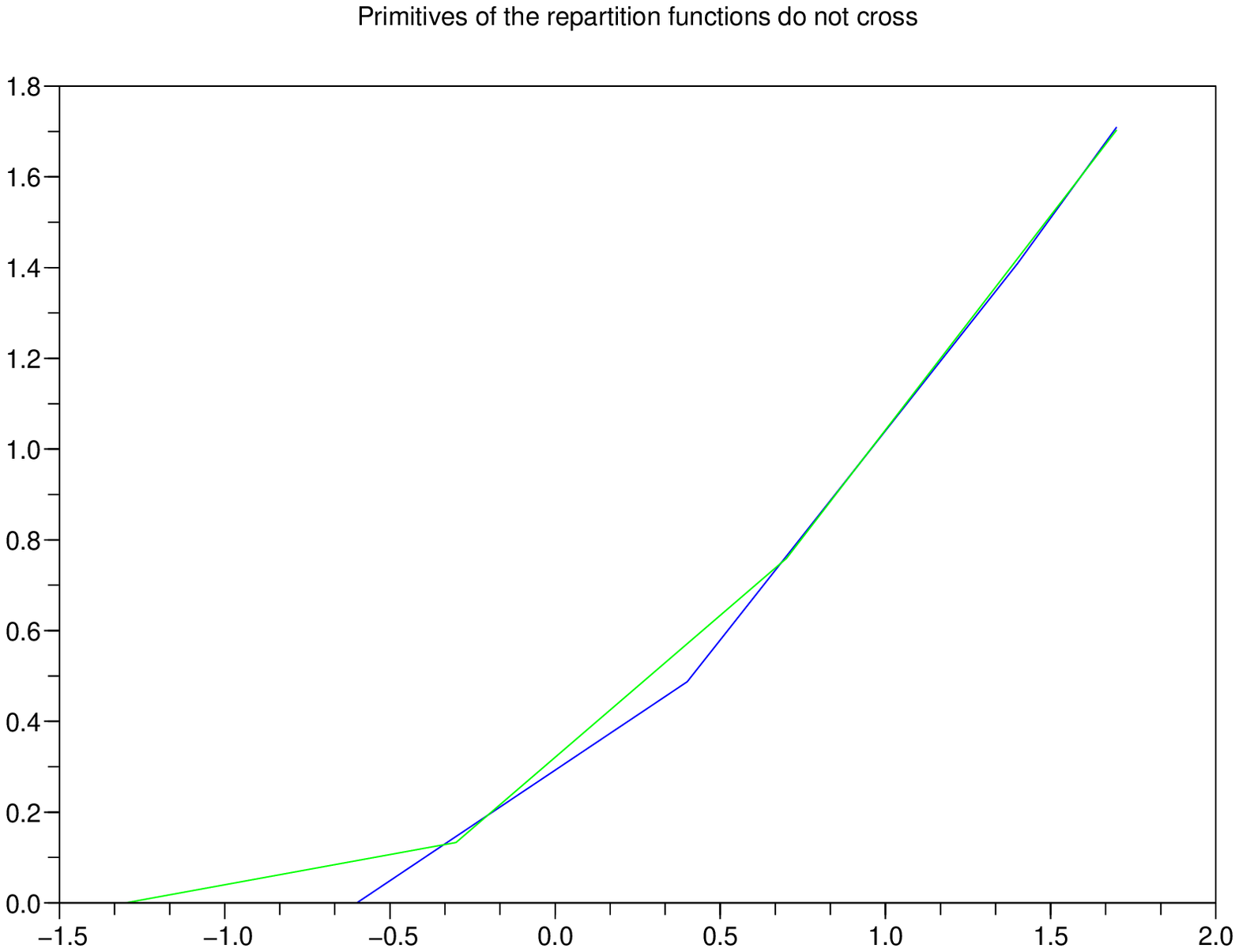}}
\end{tabular}
\caption{Illustration of the convex order: the primitives of the
repartition functions do not cross each other}
\label{fig:convex_order.eps}
\end{figure}

The convex order compares random variables according to their
``variability'' or ``spread'' \cite{Shaked-Shanthikumar:2007}.
We say that $X$ \emph{is less than} $Y$ 
\emph{in convex order}, denoted 
\begin{equation*}
X \lesser_{cx} Y \; ,
\end{equation*}
if  and only if one of the following equivalent conditions holds true:
\begin{itemize}
\item the means are equal and 
the primitive of the repartition function of $X$ is always 
below that of $Y$, that is,
$\EE(X)=\EE(Y)$ and $X \lesser_{icx} Y$,
% \item $\EE(X)=\EE(Y)$ and 
%$ \int_0^p F^{-1}(q)dq \geq  \int_0^p G^{-1}(q)dq $ for all $p\in [0,1]$, 
% where $F$ and $G$ are the repartition functions of $X$ and of $Y$, 
\item $\EE ( \varphi (X) )  \leq \EE ( \varphi (Y) )$ 
for all convex function $\varphi$.
\end{itemize}
Roughly speaking, $Y$ is more likely to take on extreme values than $X$
(see Figure~\ref{fig:convex_order.eps}).
Notice that the convex order is more demanding than the increasing
convex order since the class of ``test functions'' is 
larger: all convex functions and not only the increasing convex ones.
This is why we obtain stronger important properties that,
when $X \lesser_{cx} Y $, the means are equal $\EE(X) = \EE(Y)$,
and the variance are ordered $\vari(X) \leq \vari(Y)$.

% \begin{figure}[hptb]
% 	\centering
% \includegraphics[width=0.8\textwidth]%
% {\Figdir convex_order.eps}
% \caption{Illustration of the convex order: the primitives of the
% repartition functions do not cross each other}
% \label{fig:convex_order.eps}
% \end{figure}

\subsection{Some properties}

\begin{itemize}
\item This $icx$ and $cx$ orders are stricter than the order defined by
comparing variances:  not all pairs of random variables may be ranked.

\item Consider the class $\MM_{\mu,\sigma^2}$ of random variables having 
same mean $\mu$ and variance $\sigma^2$. Elements of
$\MM_{\mu,\sigma^2}$ cannot be compared with respect to $cx$.
Indeed, if $X \lesser_{cx} Y $ and $\vari(X)=\vari(Y)$, then $X$ and
$Y$ have the same distribution \cite[p.57]{Muller-Stoyan:2002}. 

\item Adding zero mean independent noise to a random variable increases 
variability: if $Z$ is independent of $X$ and has zero mean, then 
$X$ is less than $Y=X+Z$ in convex order.
This is a consequence of Strassen's Theorem
\cite[p.23]{Muller-Stoyan:2002}. 
More generally, without assuming independence,
$X$ is less than $Y=X+Z$ in convex order whenever
the conditional expectation \( \EE[Z|X]=0 \).

\item Consider $X$ following Normal distribution ${\cal N}(\mu,\sigma^2)$
and $Y$ following ${\cal N}(\nu,\tau^2)$.
Then, $X \lesser_{icx} Y$ if and only if 
$\mu \leq \nu$ and $\sigma^2 \leq \tau^2$,
and $X \lesser_{cx} Y$ if and only if 
$\mu = \nu$ and $\sigma^2 \leq \tau^2$
\cite[p.62]{Muller-Stoyan:2002}.

\item For $p > 0$, let us introduce 
\( \CV_p(X) \defegal \frac{\EE(X^p)^{1/p}}{\EE(X)} \)
for positive $p$-integrable random variable $X$.
For $p=2$, we have the usual \emph{coefficient of variation} 
$\CV(X) = \CV_2(X) = \sqrt{\EE(X^2)}/\EE(X)$. 
If $X \lesser_{cx} Y$, then \( \CV_2(X) \leq \CV_2(Y) \)
(in fact \( \CV_p(X) \leq \CV_p(Y) \) for all $p \geq 1$).
 
\end{itemize}

\subsection{Increasing convex order and convex order for random vectors}
\label{ss:Increasing_convex_order_and_convex_order_for_random_vectors}

We shall need to compare not only random variables but random vectors
as in \cite[p.98]{Muller-Stoyan:2002} and
\cite[p.323]{Shaked-Shanthikumar:2007}.
For this, we can no longer appeal to repartition functions.
Let $X=(X_1,\ldots,X_n)$ and $Y=(Y_1,\ldots,Y_n)$ be
random vectors with finite mean.

We say that $X$ \emph{is less than} $Y$ 
\emph{in increasing convex order}, written $X \lesser_{icx} Y$,
if and only if $\EE ( \varphi (X_1,\ldots,X_n) )  \leq 
\EE ( \varphi (Y_1,\ldots,Y_n) ) $ 
for any increasing convex function $\varphi: \RR^n \to \RR$.

We say that $X$ \emph{is less than} $Y$ 
\emph{in convex order}, written $X \lesser_{cx} Y$,
if and only if $\EE ( \varphi (X_1,\ldots,X_n) )  \leq 
\EE ( \varphi (Y_1,\ldots,Y_n) ) $ 
for any convex function $\varphi: \RR^n \to \RR$.
In this case, $X$ and $Y$ have the same mean.
\bigskip

Consider $X$ following Normal distribution ${\cal N}(\mu,\Sigma)$
and $X'$ following ${\cal N}(\mu',\Sigma')$.
Then, $X \lesser_{cx} X'$ if and only if 
$\mu = \mu'$ and $\Sigma' - \Sigma$ is non-negative definite.
The situation is not as clear cut for the $icx$ order.
If \alert{$\mu_X \geq \mu_Y$} and 
\alert{$\Sigma_X - \Sigma_Y > 0$} (non-negative definite), 
then \alert{$X \greater_{icx} Y$}.
If $X \greater_{icx} Y$, then $\mu_X \geq \mu_Y$ and 
$a^T(\Sigma_X - \Sigma_Y)a \geq 0 $ for all vector $a \geq 0$
\cite[p.100]{Muller-Stoyan:2002}.

\section{Generic results on environmental noise variability
in population dynamics matrix models}
\label{sec:Generic_results}

In what follows, we shall consider a population described at discrete
times $t=0,\ldots,\horizon$
(where \(\horizon\) is the \emph{horizon}), 
either by a scalar $\abundance(t) \in \RR$ 
or by a vector 
$\abundance(t)= \big(\abundance_1(t),\ldots,\abundance_n(t) \big) \in \RR^n$ 
which may be abundances at ages or stages.  
The \emph{population size} is 
\( \size = \| \abundance \| = \abundance_1 + \cdots + \abundance_n \).

The dynamical evolution of the population is supposed to be linear in
the sense that 
\begin{equation}
 \abundance(t+1) = \matrix\big( \varepsilon(t) \big) \abundance(t) \; ,
\quad t=0,\ldots,\horizon -1 \; ,
\label{eq:dynamics}
\end{equation}
where the matrix $\matrix$ is independent of $\abundance(t)$
(no density-dependence effect, this is why we
label such model of linear).
On the other hand, the components  $\matrix_{ij}$ 
of the matrix $\matrix$ may depend on the 
\emph{environmental factors}, a vector $\varepsilon(t)=\big(
\varepsilon_1(t), \ldots, \varepsilon_p(t) \big) \in \RR^p $ at time $t$. 

For instance, the components of the matrix $\matrix$ may depend 
linearly on the environmental factors, as in the expression
\begin{equation}
 \matrix\big( \varepsilon \big) = \left( \begin{array}{ccc}
 \overline{\matrix}_{11} + \varepsilon_{11} & \cdots &
 \overline{\matrix}_{1n} + \varepsilon_{1n} \\
\cdots & \cdots & \cdots \\
 \overline{\matrix}_{n1} + \varepsilon_{n1} & \cdots &
 \overline{\matrix}_{nn} + \varepsilon_{nn} 
\end{array} \right) 
\label{eq:linear_dependence}
\end{equation}
%This form corresponds to 
or may depend exponentially as in 
\begin{equation}
 \matrix \big( \varepsilon \big)= \left( \begin{array}{ccc}
\exp\big(\overline{\matrix}_{11} + \varepsilon_{11}\big) & \cdots &
\exp\big(\overline{\matrix}_{1n} + \varepsilon_{1n}\big) \\
\cdots & \cdots & \cdots \\
\exp\big(\overline{\matrix}_{n1} + \varepsilon_{n1}\big) & \cdots &
\exp\big(\overline{\matrix}_{nn} + \varepsilon_{nn}\big)
\end{array} \right) \; .
\end{equation}
In this latter case, the components of the matrix $\matrix$ are 
log-convex function of the environmental factors.
Recall that $f$ is a log-convex function if $f>0$ and $\log f$ is
convex. Otherwise stated, $f$ is the exponential of a convex function
(as a consequence, a log-convex function is also convex).

In \cite{Ives-Hughes:2002}, different non linear models are recalled.
When $\alpha=0$, they are matrix models without density-dependency.
Model (2a) exhibits components which are exponential in the
environmental factor, 
while they are linear in models (2c) and (2d).
Calculation shows that model (2b) has matrix components which are 
log-convex functions of the environmental factor.

We shall coin \emph{environmental scenario} a temporal sequence 
\( \varepsilon(\cdot) = 
\big( \varepsilon(0),\ldots,\varepsilon(\horizon-1) \big) \)
of environmental factors.

\begin{proposition}
Consider two environmental scenarii, one being 
\alert{more variable in increasing convex order} than the other:
\(
\big( \varepsilon^M(0),\ldots, \varepsilon^M(\horizon-1) \big) 
\greater_{icx} 
\big( \varepsilon^L(0),\ldots, \varepsilon^L(\horizon-1) \big) \).
Denote by 
\( \size^M(\horizon) = \| \matrix\big( \varepsilon^M(\horizon-1) \big) \cdots 
\matrix\big( \varepsilon^M(0) \big) \abundance(0) \| \)
and
\( \size^L(\horizon) = \| \matrix\big( \varepsilon^L(\horizon-1) \big) \cdots 
\matrix\big( \varepsilon^L(0) \big) \abundance(0) \| \)
the corresponding populations sizes. 

Assume that the components $\matrix_{ij}(\varepsilon_1,\ldots,\varepsilon_p) $ 
of the matrix $\matrix$ in~\eqref{eq:dynamics}
are nonnegative combinations %sums 
of  log-convex functions of the environmental factor 
$(\varepsilon_1,\ldots,\varepsilon_p)$.
Then, the more variable the scenario,
the more variable the population size in the sense that 
% then the corresponding population size \( \size^M(\horizon) \)
% is {more variable} in increasing convex order sense than 
\begin{equation}
\size^M(\horizon) \greater_{icx} \size^L(\horizon) 
\mtext{ and }
\log \size^M(\horizon) \greater_{icx} \log \size^L(\horizon) \; .
\end{equation}
As a consequence, \( \EE \left( \size^M(\horizon) \right) \geq
\EE \left( \size^L(\horizon) \right) \) and
\( \EE \left( \log \size^M(\horizon) \right) \geq
\EE \left( \log \size^L(\horizon) \right) \).
\end{proposition}

In a sense, 
\emph{environmental variability increases both 
mean population size and mean log-population size 
and makes them more variable}.
\medskip

\begin{proof}
The components of the vector $\abundance(\horizon)=
\matrix\big( \varepsilon(\horizon-1) \big) \cdots 
\matrix\big( \varepsilon(0) \big) \abundance(0) $ are sums of products of
nonnegative combinations of 
log-convex functions of the environmental scenario.
Therefore, by a property of log-convex functions \cite{J-Cohen:1980},
the components of the vector $\abundance(\horizon)$ are also 
log-convex functions of the environmental scenario,
and so is the population size. 
Thus, the logarithm \( \log \size(\horizon) \)
of the population size is convex
in \( \big( \varepsilon(0),\ldots,\varepsilon(\horizon-1) \big) \).
For any increasing convex function $ \varphi : \RR \to \RR$,
\( \varphi \big( \log \size(\horizon) \big) \) is convex
in \( \big( \varepsilon(0),\ldots,\varepsilon(\horizon-1) \big) \)
since convexity is preserved by left-composition with an increasing 
convex function.
We end up by using the definition of
increasing convex order for random vectors
in~\S\ref{ss:Increasing_convex_order_and_convex_order_for_random_vectors}: 
\( \EE \big[ \varphi \big( \log \size^M(\horizon) \big) \big]
\geq \EE \big[ \varphi \big( \log \size^L(\horizon) \big) \big] \).
This precisely means that 
\( \log \size^M(\horizon) \greater_{icx} \log \size^L(\horizon) \). 

Since a log-convex function is also convex,
the population size $\abundance(\horizon)$ is a sum 
of convex functions of the environmental scenario.
Then, the proof follows as above.

At last, we use the property that $X \greater_{icx} Y \Rightarrow
\EE(X) \geq \EE(Y)$ to compare the means.

\end{proof}

\begin{figure}[hptb]
	\centering
\includegraphics[width=0.8\textwidth]%
{\Figdir 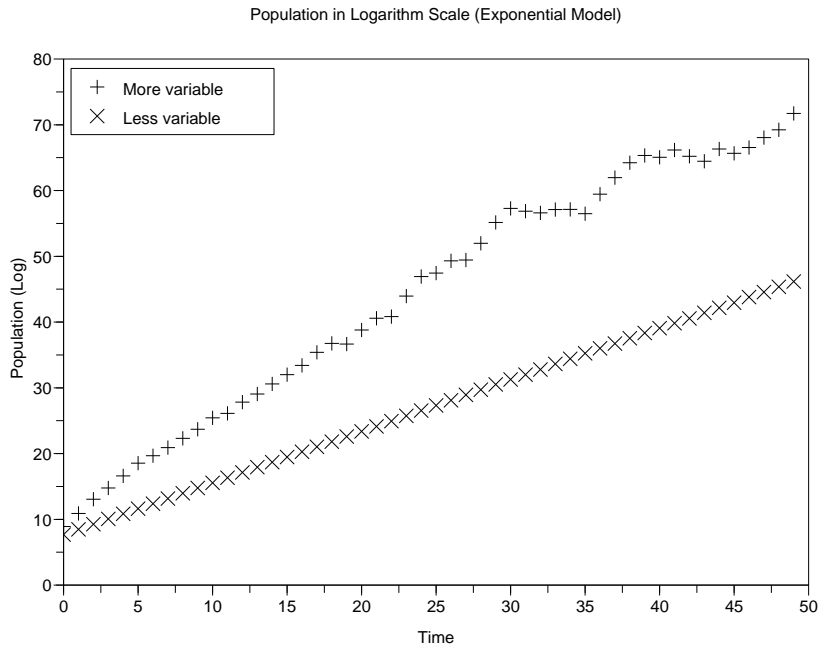}
\caption{Environmental variability increases 
mean population size for a population
  model where growth rate depends exponentialy on environmental factor}
\label{fig:exp_variability.eps}
\end{figure}

Instead of total population, the result would still hold true with any 
positive weighted combination
$a_1 \abundance_1 +\cdots + a_k \abundance_k $ where $a_i \geq 0$,
or with 
$\log(a_1 \abundance_1 +\cdots + a_k \abundance_k )$ where $a_i \geq 0$,
\bigskip

As an illustration, 
consider the following scalar dynamic equation for population size
\(
 \abundance(t+1) = \exp \big(r + \varepsilon(t) \big)
\abundance(t) \; ,
\)
for which we have 
\(
 \abundance(\horizon) = 
\exp \big( r\horizon + \varepsilon(0) + \cdots + \varepsilon(\horizon-1) \big) 
\abundance(0) \; .
\)
Hence, both \( \abundance(\horizon) \) and 
\( \log \abundance(\horizon) \) are convex 
functions of the environmental scenario
$\varepsilon(\cdot)=(\varepsilon(0),\ldots,\varepsilon(\horizon-1))$, 
so that environmental variability increases mean population size 
% and mean logarithmic population size 
as may be seen in Figure~\ref{fig:exp_variability.eps}.
Indeed, the mean population size generated by a more variable environment
is above the one by a less variable environment, for all times.

\section{Conclusion}
\label{sec:Conclusion}

We have used another notion of \alert{variability} than 
the widely used variance or coefficient of variation, namely the
so-called convex orders.
We think that such partial orders may be of interest in theoretical
ecology beyond this specific application. 

To compare our approach with the literature, notice that,
though we consider matrix population models, we make 
no ergodic assumption on the stochastic process 
$\matrix_0$, $\matrix_1$, \ldots
However, we make separate assumptions, on the one hand
on the environmental factors 
\alert{$\varepsilon(0),\ldots, \varepsilon(t)$} 
and, on the other hand, on the functional dependence 
\alert{$\matrix_t = \matrix\big( \varepsilon(t) \big) $}.

With this approach, we obtain generic results which 
are \alert{not asymptotic in time}, 
but \alert{valid at any time $t$} and for a large class of 
functional dependence on the uncertainties.

Though the approach is different, our conclusions are consistent
with the cases presented
in Sect.~\ref{sec:Convex_orders}.
We extend the observation of Lande that,
when adding environmental noise to growth rate on the log scale,
\alert{environmental  stochasticity increases the mean multiplicative
  growth rate} to matrix models.
As to Tuljapurkar's asympotic approximation,
we arrive at a different conclusion because his assumptions
correspond to a matrix $\matrix$ depending
linearly on the environmental factors as in~\eqref{eq:linear_dependence},
and our result does not cover this case.

Our general conclusion is, therefore, that the analytical dependence on 
environmental factors cannot be overlooked when trying to tackle the
influence of variability.
However, as shown in this paper,
specific analytical dependence coupled with appropriate notion
of  variability lead to wide generic results, valid for all times and
not only asymptotically, and requiring no assumptions of
stationarity, of normality, of independency, etc.

\paragraph{Acknowledgements.}
Some years ago, Shripad Tuljapurkar encouraged me to develop the general
ideas I exposed to him after one of his talks in Paris, and I thank him
for this.
I want to thank Michel Loreau and Claire de Mazancourt for welcoming me 
at the Dept of Biology, Mc Gill, Montreal, Canada.
Fruitful discussions with them and with the participants to a seminar in
August 2008 helped me shape my ideas.
I also want to thank Tim Coulson and the participants to the 
Ecology and Evolution Seminar Series, Silwood Park campus, United
Kingdom on March 2009.

\newcommand{\noopsort}[1]{} \ifx\undefined\allcaps\def\allcaps#1{#1}\fi


\begin{thebibliography}{1}

\bibitem{Lande-Engen-Saether:2003}
R.~Lande, S.~Engen, and B.-E. Saether.
\newblock {\em Stochastic population dynamics in ecology and conservation}.
\newblock Oxford series in ecology and evolution, 2003.

\bibitem{Tuljapurkar:1990}
S.~Tuljapurkar.
\newblock {\em Population Dynamics in Variable Environments}.
\newblock Springer-Verlag, Berlin, 1990.
\newblock Lecture Notes in Biomathematics.

\bibitem{Caswell:2001}
H.~Caswell.
\newblock {\em Matrix Population Models}.
\newblock Sinauer Associates, Sunderland, Massachussetts, second edition, 2001.

\bibitem{Muller-Stoyan:2002}
A.~Muller and D.~Stoyan.
\newblock {\em Comparison Methods for Stochastic Models and Risk}.
\newblock John Wiley and Sons, New York, 2002.

\bibitem{Shaked-Shanthikumar:2007}
Moshe Shaked and J.~George Shanthikumar.
\newblock {\em Stochastic Orders}.
\newblock Springer-Verlag, Berlin, 2007.

\bibitem{Ives-Hughes:2002}
Anthony~R. Ives and Jennifer~B. Hughes.
\newblock General relationships between species diversity and stability in
  competitive systems.
\newblock {\em Am. Nat.}, 159(4):388--395, April 2002.

\bibitem{J-Cohen:1980}
Joel~E. Cohen.
\newblock Convexity properties of products of random nonnegative matrices.
\newblock {\em Proc. Nat. Acad. Sci. USA}, 77:3749--3752, 1980.

\end{thebibliography}
\end{document}